\documentclass[smallextended]{svjour3}
\newcommand{\be}{\begin{equation}}
\newcommand{\ee}{\end{equation}}
\newcommand{\beq}{\begin{eqnarray}}
\newcommand{\eeq}{\end{eqnarray}}
\newcommand{\nbeq}{\begin{eqnarray*}}
\newcommand{\neeq}{\end{eqnarray*}}

\begin{document}

\title{Characterization of exponential distribution through bivariate regression of record values revisited   }
\author{George P. Yanev }
\institute{George P. Yanev \at
             School of Mathematical and Statistical Sciences, The University of Texas Rio Grande Valley\\
              1201 W. University Drive, Edinburg, Texas, 78539 USA\\
              Tel.: (956) 665-3632, \email{george.yanev@utrgv.edu}           }

\maketitle

\begin{abstract}
It is shown that the exponential is the only distribution which satisfies a
certain regression equation. This characterization equation involves the conditional expectation (regression function) of a
record value given a pair of record values, one previous and one future, as covariates. The underlying distribution is
exponential if and only if the above regression equals the expected value of an appropriately
defined Beta distributed random variable. In a particular case, the expected value of the Beta variable reduces to a 
weighted average of the covariates.
\end{abstract}


\titlerunning{Characterization through regression on two covariates} % if too long for running head

\institute{George P. Yanev \at
              Department of Mathematics, The University of Texas Rio Grande Valley\\
              1201 W. University Drive, Edinburg, Texas, 78539 USA\\\bar{}
              Tel.: (956) 665-3632, \email{george.yanev@utrgv.edu}           }

\section{Introduction}
I first met Professor Ahsanullah a.k.a. Moe at a conference in 2002. Later, he spent the academic year 2005-2006 
at the Department of mathematics of the University of South Florida as a Visiting Professor. 
Despite some health problems, Moe taught  two classes, did research, 
and advised a Ph.D. student, who later defended successfully his thesis. At one of the
weekly meetings of the Probability and Statistics seminar, Moe posed the open question of 
characterizing probability distributions by bivariate regression of record values. Under his leadership, 
M. Beg (also visiting) and myself started working on that problem. The results of this collaboration appeared in 
Yanev et al. (2008) and were extended in Yanev~(2012). Here, we revisit the bivariate regression problem, obtaining an alternative 
form for the right-hand side of the  characterization equation and providing some additional insight.  

To formulate and discuss the obtained results, we need to introduce the following  notation. 
Let $X_1, X_2, \ldots$ be independent copies of a random variable 
$X$ with an absolutely continuous distribution function $F(x)$. One observation in
a discrete time series is called a (upper) record value if it exceeds all
previous observations, i.e., $X_j$ is a (upper) record value if $X_j>X_i$
for all $i<j$. More precisely, let us define the classical (upper) record times $T_n$ and (upper) record values $R_n$ as follows:
$T_1=1$, $R_1=X_1$, and then $T_{n+1}=\min \{j:X_j>R_n\}$, $R_{n+1}=X_{T_{n+1}}$, $n=1,2\ldots$ (see Ahsanullah and
Nevzorov (2015), p.46).

Let $F(x)$ be the exponential cumulative distribution function
\begin{equation}  \label{exp_type}
F(x)=1-e^{\displaystyle -c(x-l_F)} \qquad (-\infty<l_F\le x),
\end{equation}
where $c>0$ is an arbitrary constant.
The distribution (\ref{exp_type})
with $l_F>0$  appears, for example, in reliability studies where $l_F$
represents the guarantee time; that is,  failure cannot occur before $l_F$
units of time have elapsed (see Barlow and Proschan (1996), p.13).

We study characterizations of exponential distributions in terms of the
regression of one record value with two other record values as covariates.
More precisely,  we examine the regression function for $1\le s\le n-1$ and $r\ge 1$ 
\be \label{regr}
E[\psi(R_n)\ |\ R_{n-s}=u, R_{n+r}=v ] \qquad (l_F\le u<v),
\ee
where $\psi$ is a function that satisfies certain regularity conditions. For a connection with the doubly-censored regression $E[\psi(X)\ |\ x<X\le y]$, we refer to Pakes (2004), Section 5.

Let us introduce 
a four-parameter (generalized) Beta random variable $B_{r,s}(u,v)$ with probability density function (see Johnson et al. (1995) vol. 2, p.210) for $r>0$ and $s>0$
\be \label{beta}
f_B(y)=\frac{1}{B(r,s)}\frac{(y-u)^{s-1}(v-y)^{r-1}}{(v-u)^{r+s-1}} \qquad (u\le y\le v),
\ee
where $B(\cdot,\cdot)$ is the Beta function. The family of distributions (\ref{beta}) includes the uniform ($s=r=1$) and the power function ($s=1$ or $r=1$) distributions as special cases. 

The following result should be known, however I could not find it formulated anywhere.

\vspace{0.3cm} {\bf Proposition}\  If $X$ is exponential with (\ref{exp_type}) and $E[\psi(B_{r,s}(u,v))]<\infty$,  then for  $r\ge 1$, $1 \le s\le n-1$, and $0<u<v$
\be \label{regression1}
E\left[\psi(R_n)\ |\ R_{n-s}=u, R_{n+r}=v\right] = E\left[\psi(B_{r,s}(u,v))\right].
\ee
{\bf Remark}. Define the increments of record values  as $\Delta_{m,n}:=R_n-R_m$ for $n>m$.
%\[
%\Delta_{n-s,n}:=R_n-R_{n-s}\qquad \mbox{and}\qquad \Delta_{n,n+r}:=R_{n+r}-R_n.
%\]
It is known (see Ahsanullah and Nevzorov (2015), p.65) that, if the underlying distribution is exponential, then for any $k$,
$
R_k \stackrel{d}{=} E_1+E_2+\ldots +E_k,
$
where $E_i$, $i=1,2,\ldots, k$ are independent and unit exponential random variables. Therefore, $\Delta_{n-s,n}$ and $\Delta_{n,n+r}$ are independent and distributed $Ga(s,1)$ and $Ga(r,1)$, respectively,
%\[
%\Delta_{n-s,n}=E_{n-s+1}+\ldots+E_n \sim Ga(s,1)\quad \mbox{and}\quad 
%\Delta_{n,n+r}=E_{n+1}+\ldots+E_{n+r}\sim Ga(r,1),
%\]
where $Ga(\cdot,\cdot)$ denotes the Gamma distribution. Thus, using a well-known property of Beta distribution 
(see Johnson et al. (1994), vol.1, p.350) the equation (\ref{regression1}) can be written as
\[
E\left[\psi(R_n)\ | \ R_{n-s}=u, R_{n+r}=v\right]= E\left[ \psi\left(\frac{u\Delta_{n,n+r}+v\Delta_{n-s,n}}
{\Delta_{n-s,n+r}} \right)   \right].
\]
Next, we shall address the question if, under some regularity assumptions on $F$ and $\psi$ and their derivatives, 
(\ref{regression1}) is also a sufficient condition for (\ref{exp_type}).
Bairamov et al. (2005) consider  (\ref{regr}) in the particular case  when both covariates are adjacent (one spacing away) to $R_n$. 
(See also Bairamov and Ozkal (2007) for similar result about order statistics.) 
They prove, under some regularity conditions, that $F$
is exponential if and only if for a function $\psi$
\be \label{BAP}
E\left[\psi(R_n)\ |\ R_{n-1}=u, R_{n+1}=v \right]=\frac{1}{v-u}\int_u^v \psi(t)\, dt \qquad (l_F\le u<v).
\ee
Observing that the right-hand side of (\ref{BAP}) equals $E[\psi(B_{1,1}(u,v))]$, 
one can rewrite (\ref{BAP}) as
\be \label{newBAP}
E\left[\psi(R_n)\ |\ R_{n-1}=u, R_{n+1}=v \right]=E[\psi\left(B_{1,1}(u,v)\right)]\qquad (l_F\le u<v).
\ee
Further on we will assume that the function $g(x)$ satisfies the following conditions for some positive integers $r$ and $s$:
\begin{description}
  \item[(i)] $g^{(r+s)}(x)$ exists and is continuous in $(l_F, \infty)$;
  \item[(ii)] $|g^{(r+s)}(l_F+)|<\infty$ for $r=2$;   $|g^{(r+s+r-1)}(l_F+)|<\infty$ for $r\ge 3$;
  \item[(iii)] $g^{(r+s)}(l_F+)\ne 0$.
\end{description}
Denote the cumulative hazard function of $X$ by $H(x):=-\ln(1-F(x))$ for $x\ge l_F$. Let $F(x)$ satisfies the following conditions for some positive integers $r$ and $s$:
\begin{description}
  \item[(iv)] $F^{(m)}(x)$ for $m=\max\{s,r\}$ exists and is continuous in $(l_F, \infty)$;
  \item[(v)] $H'(x)$ is nowhere constant in a small interval $(l_F, l_F+\varepsilon)$ for $\varepsilon>0$;
  \item[(vi)] $H'(l_F+)>0$ and $|H^{(m)}(l_F+)<\infty$ for $m\le \max \{3,r\}$.
\end{description}
Extending (\ref{newBAP}) to covariate record values  non-adjacent to $R_n$, we obtain our main result.

\vspace{0.3cm} {\bf Theorem}\ If (i)-(vi) hold and for $1\le s\le n-1$, $r\ge 1$, and $l_F<u<v$
\be  \label{main}
\hspace{-0.5cm}E\left[g^{(r+s-1)}(R_n)\ |\ R_{n-s}=u, R_{n+r}=v\right] =E\left[g^{(r+s-1)}(B_{r,s}(u,v))\right],
\ee
then $X$ is exponential with (\ref{exp_type}) for some $c>0$.

Setting $g(x)=x^{r+s}/(r+s)!$, hence 
$g^{(r+s-1)}(x)=x$, and taking into account that $E[B_{r,s}(u,v)]=(ur+vs)/(r+s)$,
one can see that the above results imply the following.

\vspace{0.3cm} {\bf Corollary }\ Let $n$, $r$, and $s$ be integers, such that $%
1\le s\le n-1$ and $r\ge 1$. Suppose that the assumptions (iv)-(vi) of the Theorem
hold. Then $F(x)$ is exponential (\ref{exp_type}) with $
c=H'(l_F+)$
if and only if
\begin{equation}  \label{thmg}
E[R_n|R_{n-s}=u, R_{n+r}=v]=\frac{r}{r+s}u+\frac{s}{r+s}v \qquad (l_F\le u<v).
\end{equation}
Note that the right-hand side of (\ref{thmg}) is a linear function of $u$ and $v$. It is a weighted average of the two covariate record values, 
where the weight of each covariate is proportional to the distance, in number of spacings,  from $R_n$ to the 
other covariate. In particular, for any $r$, such that $1\le r\le n-1$, (\ref{thmg}) simplifies to
\[
E[R_n|R_{n-r}=u, R_{n+r}=v]=\frac{u+v}{2} \qquad (l_F\le u<v).
\]
In Sections 2 and 3, we shall prove the Proposition and the Theorem, respectively. The last section includes some concluding remarks.

\section{Proof of the Proposition}\
Using the Markov dependence of record values (e.g., Nevzorov (2001), p.68),
for the conditional density of $R_n$
given $R_{n-s}=u$  and $R_{n+r}=v$ for $u\le t\le v$ 
we obtain
\beq \label{cond_density}
f_{n|n-s,n+r}(x|u,v) & = & \frac{f_{n+r|n-s,n}(v|u,x)f_{n-s,n}(u,x)}{f_{n-s,n+r}(u,v)}\\
    & = & \frac{f_{n+r|n}(v|x)f_{n-s,n}(u,x)}{f_{n-s,n+r}(u,v)} \nonumber \\
    & = & \frac{f_{n,n+r}(x,v)f_{n-s,n}(u,x)}{f_n(x)f_{n-s,n+r}(u,v)}. \nonumber
\eeq
Assuming (\ref{exp_type}), we have (e.g., Ahsanullah and Nevzorov (2015), p.80) 
\be \label{joint_density}
f_n(x)=\frac{c^n(x-l_F)^n}{n!}f(x),\qquad  f_{m,n}(x_m, x_n)=\frac{c^n(x_m-l_F)^m(x_n-x_m)^{n-m-1}}{m!(n-m-1)!}f(x_n).         
 \ee
Combining (\ref{cond_density}) and (\ref{joint_density}), we obtain
\[
f_{n|n-s,n+r}(t|u,v)= \frac{1}{B(r,s)}\frac{(t-u)^{s-1}(v-t)^{r-1}}{(v-u)^{r+s-1}},
\]
which is the probability density function of a four-parameter Beta distribution.
Therefore, 
\nbeq
E\left[\varphi(R_n)\ |\ R_{n-s}=u, R_{n+r}=v\right] & = &
    \frac{1}{B(r,s)}\int_u^v \varphi(t)\frac{(t-u)^{s-1}(v-t)^{r-1}}{(v-u)^{r+s-1}}\, dt \\
    & = &
 E\left[\varphi(B_{r,s}(u,v))\right],
\neeq
which proves (\ref{regression1}).

\section{Proof of the Theorem}
For a given function $g(x)$ and non-negative integers $i$
and $j$, define for $u\ne v$
\[  
\ _iM_j(u,v):=\frac{\partial^{i+j}}{\partial
u^i\partial v^j}\left(\frac{g(v)-g(u)}{v-u} \right).            
\]
Under the assumptions of the Theorem, it was proven in Yanev~(2012) that if
\be \label{2012}
E\left[g^{(r+s-1)}(R_n)\ |\ R_{n-s}=u, R_{n+r}=v\right]=\frac{1}{B(r,s)}\ \  _{r-1}M_{s-1}(u,v),
\ee
then  $X$ is exponential with (\ref{exp_type}). Therefore, to prove the Theorem, it is sufficient to show that the 
right-hand sides of (\ref{main}) and (\ref{2012}) are equal, i.e., for $r\ge 1$ and $1\le s\le n-1$
\be \label{exp_M}
E[g^{(r+s-1)}(B_{r,s}(u,v))]=\frac{1}{B(r,s)}\ \  _{r-1}M_{s-1}(u,v).
\ee
It is not difficult to verify (\ref{exp_M}) for $r=s=1$. Indeed, referring to (\ref{beta}), we have
\[
E[g'(B_{1,1}(u,v))]=\frac{1}{B(1,1)}\int_u^v \frac{g'(t)}{v-u}\, dt=\frac{1}{B(1,1)}\frac{g(v)-g(u)}{v-u}.
\]
Next, assuming (\ref{exp_M}) for $r=1$ and $1\le s=i\le n-2$, we shall prove it for $r=1$ and $s=i+1$. 
One can verify (see Lemma~1 in Yanev et al. (2008)) the following identity between the derivatives $g^{(j)}(x)$ and $\ _0M_j(u,v)$ 
for $j\ge 1$ and $u<v$
\begin{equation}  \label{lemma_old}
g^{(j)}(v)=(v-u)\ _0M_j(u,v)+j\ _0M_{j-1}(u,v).  
\end{equation}
Using the induction assumption
\[
E\left[g^{(i)}(B_{1,i}(u,v))\right]=i\ _0M_{i-1}(u,v),
\]
we obtain
\nbeq
E[g^{(i+1)}(B_{1,i+1}(u,v))] & = & \frac{1}{B(1,i+1)}\int_u^v g^{(i+1)}(t)\frac{(t-u)^i}{(v-u)^{i+1}}\, dt \\
    & = & \frac{1}{B(1,i+1)}\frac{1}{(v-u)^{i+1}}\left[ g^{(i)}(v)(v-u)^i-i\int_u^v g^{(i)}(t)(t-u)^{i-1}\, dt\right] \\
     & = & \frac{1}{B(1,i+1)}\frac{1}{(v-u)^{i+1}}\left[ g^{(i)}(v)(v-u)^i-(v-u)^i 
     E[g^{(i)}(B_{1,i}(u,v))]\right] \\
    & = & \frac{1}{B(1,i+1)}\frac{1}{v-u}\left[ g^{(i)}(v)-i\ _0M_{i-1}(u,v)\right] \\
    & = & \frac{1}{B(1,i+1)}\ _0M_i(u,v),
 \neeq
where the last equality follows from (\ref{lemma_old}). This proves (\ref{exp_M})
for $r=1$ and any $1\le s\le n-1$. Similarly, one can prove (\ref{exp_M}) for $s=1$ and
any $r > 1$, i.e.,
\[
E\left[g^{(r)}(B_{r,1}(u,v))\right]=r\ _{r-1}M_0(u,v).
\]
To complete the proof of (\ref{exp_M}), it remains to establish it for any $r\ge 2$ and $2\le s\le n-1$. Assuming (\ref{exp_M})
for $r=j$ and any fixed $2\le s\le n-1$, we shall prove it for $r=j+1$ and $2\le s\le n-1$, i.e., we shall prove that
\[
E[g^{(s+j)}(B_{j+1,s}(u,v))]=\frac{1}{B(j+1,s)}\ _jM_{s-1}(u,v).
\]
Since
\[
E[g^{(s+j)}(B_{j+1,s}(u,v))]= \frac{1}{B(j+1,s)}\int_u^v g^{(s+j)}(t)\frac{(t-u)^{s-1}(v-t)^j}{(v-u)^{s+j}}\, dt,
\]
it is sufficient to prove that
\be \label{integral}
I(j+1,s):=\int_u^v g^{(s+j)}(t)(t-u)^{s-1}(v-t)^j\, dt=(v-u)^{s+j}\ _jM_{s-1}(u,v),
\ee
provided that (induction assumption)
\be \label{integral1}
I(j,s):=\int_u^v g^{(s+j-1)}(t)(t-u)^{s-1}(v-t)^{j-1}\, dt=(v-u)^{s+j-1}\ _{j-1}M_{s-1}(u,v).
\ee
Integrating (\ref{integral1}) by parts, we obtain
\nbeq
I(j,s)=\frac{1}{j}\int_u^v g^{(s+j)}(t)(t-u)^{s-1}(v-t)^j\, dt + \frac{s-1}{j}\int_u^v g^{s+j-1}(t)(t-u)^{s-2}(v-t)^j\, dt.
\neeq
Hence,
\[ 
I(j+1,s)=jI(j,s)-(s-1)I(j+1,s-1).
\]
Iterating last equation, we have
\beq \label{integral2}
I(j+1,s) & = & jI(k,j)-j(k-1)I(k-1,j)+j(k-1)(k-2)I(k-2,j+1)\\
	& = & j\sum_{i=0}^{s-2}a_iI(j,s-i)+a_{s-1}\int_u^v g^{(j+1)}(t)(v-t)^j\, dt, \nonumber
\eeq
where $a_i=(-1)^i(s-1)!/(s-1-i)!$ for $i=1,2,\ldots, s-1$; $a_0=1$. Observe that (see Lemma~1 in Yanev et al. (2008)) for $i, j\ge 1$ and $v>u$
\be \label{iden2}
i\ _{i-1}M_j(u,v)=(v-u)\ _iM_j(u,v)+j\ _iM_{j-1}(u,v).
\ee
Finally, from (\ref{integral2}), using repeatedly (\ref{iden2}), we obtain
\nbeq
\frac{I(j+1,s)}{(v-u)^{j+1}} & = & j\sum_{i=0}^{s-2}a_i(v-u)^{s-2-i}\ _{j-1}M_{s-1-i}(u,v)+a_{s-1}\ _jM_0(u,v)\\
    & = & j\sum_{i=0}^{s-3}a_i(v-u)^{s-2-i}\ _{j-1}M_{s-1-i}(u,v)+a_{s-2}[j\ _{j-1}M_1(u,v)-\ _jM_0(u,v)] \\
    & = & j\sum_{i=0}^{s-3}a_i(v-u)^{s-2-i}\ _{j-1}M_{s-1-i}(u,v)+a_{s-2}(v-u)\ _jM_1(u,v) \\
    & & \cdots  \\
    & = & j(v-u)^{s-2}\ _{j-1}M_{s-1}(u,v)-(s-1)(v-u)^{s-2}\ _jM_{s-2}(u,v)\\
    & = & (v-u)^{s-1}\ _jM_{s-1}(u,v).
\neeq
This implies (\ref{integral}) and thus proves (\ref{exp_M}) for any $r\ge 2$ and $2\le s\le n-1$. The proof of the theorem is complete.

\section{Concluding remarks}\
The main result in this paper is a characterization of the exponential distribution via a
bivariate regression relation of record values. Introducing an appropriate, generalized Beta distributed, random variable, 
we simplify the characteristic equation obtained  previously by Yanev~(2012). 

The regularity assumptions on the functions $F$ and $\psi$ and their derivatives
in the Theorem are the same as those in Yanev~(2012). Some of these conditions are quite technical and are needed to 
reach a contradiction in Yanev's~(2012) proof.
Using a different technique of proof, for example utilizing differential equations as in Bhatt (2013) or general integral equations, one might be able to weaken these assumptions.
Another question of interest is whether the presented characterization results can be extended to 
regression relations of order variables from other sub-classes of the generalized order statistics.

\vspace{0.2cm}{\bf Acknowledgment} 
I thank the anonymous referee for the useful suggestions, which improve the presentation. 
The author was partially supported by the NFSR
190 at the MES of Bulgaria, Grant No DFNI-I02/17 while on leave from the
Institute of Mathematics and Informatics at the Bulgarian Academy of Sciences.

\end{document}